\documentclass[11pt]{article}

\usepackage[utf8]{inputenc}
\usepackage[T1]{fontenc}
\usepackage{times} 
\usepackage{graphicx}
\usepackage{subcaption}
\usepackage{xcolor}
\usepackage{amssymb,amsmath,amsthm}
\usepackage{placeins}
\usepackage{hyperref}
\usepackage{pdfpages}
\usepackage{authblk}

\usepackage[margin=1in]{geometry}

\title{Optimal Control Strategies for Multi-Agent Sheep Herding}

\author[1,2]{Drake Brown}
\author[1,3]{Trevor Garrity}
\author[1,4]{Daniel Perkins}
\author[1,4]{Davis Hunter}
\author[1]{Wyatt Pochman}

\affil[1]{Brigham Young University}
\affil[2]{University of Utah}
\affil[3]{University of Maryland}
\affil[4]{Bredesen Center, University of Tennessee}

\begin{document}

\maketitle

\begin{abstract}
We develop a cost functional and state-space equations to model the problem of herding $m$ sheep to the origin using $n$ dogs. Our initial approach uses \texttt{solve\_bvp} to approximate optimal control trajectories. But this method often fails to converge due to the system’s high dimensionality and nonlinearity. However, with a well-chosen initial guess and carefully selected hyperparameters, we succeed in getting \texttt{solve\_bvp} to converge. We also explore alternatives including the shooting method and linearization with the iterative Linear Quadratic Regulator (iLQR). While the shooting method also suffers from poor convergence, the linearized iLQR approach proves more scalable and successfully handles scenarios with more agents. However, it struggles in regions where dogs and sheep are in close proximity, due to strong nonlinearities that violate the assumptions of local linearization. This leads to jagged, oscillatory paths and slow convergence, particularly when the number of sheep exceeds the number of dogs. These challenges reveal key limitations of standard numerical techniques in multi-agent control and underscore the need for more robust, nonlinear strategies for coordinating interacting agents.
\end{abstract}

\section{Introduction}

Shepherds have long considered the question, \textit{what are the optimal routes for my $m$ dogs to take to corral my $n$ indignant sheep into their pen?} Assuming dogs behave optimally in real life, this is closely related to the question, \textit{how many dogs $m$ are required to successfully corral $n$ sheep into their pen?}

Exploiting the facts that (i) sheep tend to flock together in groups and (ii) sheep are slower than dogs, shepherds have universally determined (through real-life experience) that there exists a solution $m$ such that $m\ll n$. Here we propose a mathematical approach that allows heuristic estimation of the required parameter $m$ through simulations.

While inspired by shepherding, the implications of solving this control problem extend to many other scenarios where controllable agents are tasked to guide a another group of agents toward a desired destination. For example, solutions could be naturally applied to optimal control in robotic search-and-rescue missions, oil spill containment, and charged particle manipulation.

\section{Background}

In general, previous researchers have applied heuristic algorithms to predict dog behavior. Strömbom et al. present a biologically-inspired model which attempts to replicate real-life herding behaviors, particularly the observed herding phases of ``driving" and ``collecting" \cite{Strombom2014}. Dog behavior is derived from the modeling choices that sheep are attracted to one another and repelled by dogs; the algorithm explicitly causes the dogs to cluster the sheep into a group (``collecting") before guiding the entire flock toward the target (``driving"). Results indicate the dogs often fail to corral the sheep in the case that $m\ll n$.

More recently, Lama and Di Bernardo consider the herding problem under the dual assumptions that (i) dogs are only aware of nearby sheep and (ii) the sheep do not naturally flock together \cite{Lama2024}. Sheep behavior is modeled by the same equations governing particle interaction in active matter \cite{PhysRevE.104.044613}; dog behavior is determined by imposing a fixed strategy.

Building on existing works, our goal is to explore the shepherding problem through the lens of optimal control. In particular, we determine the optimal dog strategy with respect to a chosen cost functional. By modeling the system with state-space dynamics and applying Pontryagin’s Maximum Principle, we apply boundary value problem solvers and the iterative Linear Quadratic Regulator (iLQR) to derive efficient, theoretically grounded strategies for herding.

\section{Modeling}

\subsection{Modeling Assumptions}

\begin{figure}
\begin{center}
    \includegraphics[width=0.6\textwidth]{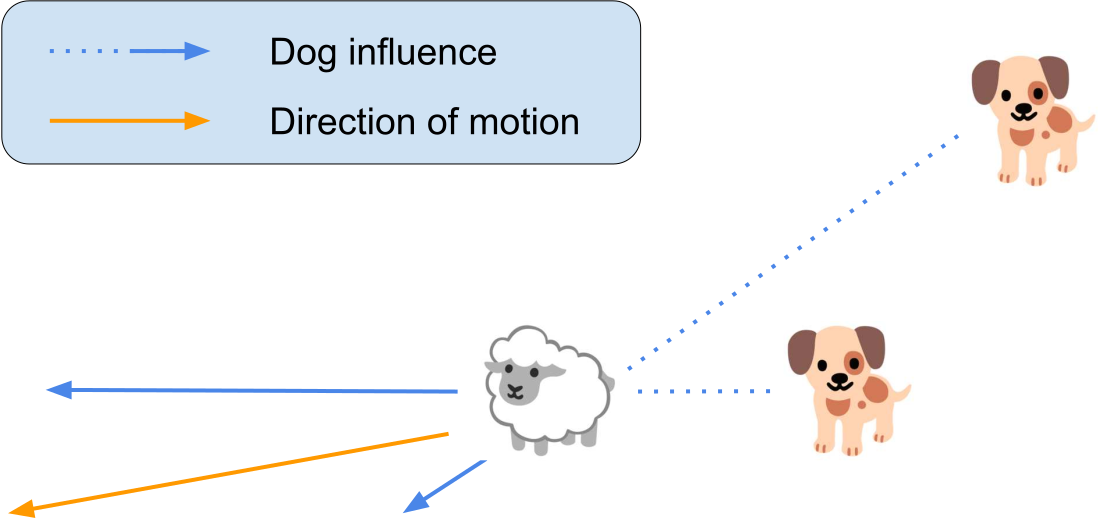}
\end{center}
     \caption[Behavior near $1/\sqrt{5}$]{A visualization of the sheep’s strategy when herded by $m=2$ dogs and $n=1$ sheep. The sheep move away from all nearby dogs, with their velocities weighted inversely by each dog’s distance—giving more influence to those that are closest.}
\label{fig:sheep_strategy}
\end{figure}

To tackle this problem, we begin by considering the behavior of the sheep. We assume each sheep wishes to maximize their distance from the dogs at each time $t$, especially moving away from the closest dog(s) (see Figure~\ref{fig:sheep_strategy}). Because sheep commonly travel in groups, we enforce no constraint on the distances between sheep, but we also do not impose any attractive forces between sheep. These modeling decisions are consistent with assumptions (i) and (ii) of \cite{Lama2024}.

Setting up the environment, we suppose the sheep and dogs are constrained to move on a two-dimensional plane, with the sheep pen at the origin. Let sheep $i$'s position be given by the two-dimensional vector $\mathbf{s}^{(i)}=(s_{x}^{(i)},s_{y}^{(i)})$, and let dog $j$'s position be given by the two-dimensional vector $\mathbf{d}^{(j)}=(d_{x}^{(j)},d_{y}^{(j)})$. We model the sheep strategy from Figure~\ref{fig:sheep_strategy} using the following second-order ODE giving the acceleration of sheep $i$: 
\begin{equation}
    (\mathbf{s}^{(i)})'' = \sum_{j=1}^m \frac{\mathbf{s}^{(i)}-\mathbf{d}^{(j)}}{(\|\mathbf{s}^{(i)}-\mathbf{d}^{(j)}\|_2^2+\epsilon)^{\lambda/2}}
\label{eq:sheep}
\end{equation}
In the case $\lambda = 3$, this acceleration is proportional to Coloumb's law governing the acceleration induced by charged particles on another like-charged particle, with the $\epsilon$ term added to avoid numerically instability that would otherwise result as the distances approach zero. Note that this $\epsilon$ term implicitly bounds the maximum possible acceleration of the sheep, with a larger value of $\epsilon$ giving a smaller maximum acceleration.

\subsection{Cost Functional and State-Evolution Equation}

The control problem is to find the optimal acceleration $\textbf{u}^{(j)}=(\mathbf{d}^{(j)})''$ of each dog $j$, given the following goals:
\begin{itemize}
\item The dogs wish to guide the sheep towards the origin (sheep pen);
\item The dogs wish to stay relatively close to the origin;
\item The dogs wish to limit their acceleration to preserve energy.
\end{itemize}

To accomplish our goals, for some fixed time $t_f$, we wish to solve the following problem:
\begin{align}
\text{minimize }&\int_0^{t_f}\left(\alpha\sum_{i=1}^n\|\mathbf{s}^{(i)}\|_2^2+\beta\sum_{j=1}^m\|\mathbf{d}^{(j)}\|_2^2+\sum_{j=1}^m\|\mathbf{u}^{(j)}\|_2^2\right)dt \label{eq:Cost Functional} \\
\text{subject to }&\begin{pmatrix} d_x^{(j)} \\ d_y^{(j)} \\ (d_x^{(j)})' \\ (d_y^{(j)})' \end{pmatrix}'=\begin{pmatrix} (d_x^{(j)})' \\ (d_y^{(j)})' \\ u_x^{(j)} \\ u_y^{(j)} \end{pmatrix}\quad \forall 1\leq j\leq m,\qquad&\label{eq:dog_evolution}
\end{align}
\begin{align}
&\begin{pmatrix} d_x^{(j)}(0) \\ d_y^{(j)}(0) \end{pmatrix} = \mathbf{d}_{0}^{(j)} \quad \forall 1\leq j\leq m,\qquad& \\ &\begin{pmatrix} s_{x}^{(i)} \\ s_y^{(i)} \\ (s_{x}^{(i)})' \\ (s_{y}^{(i)})' \end{pmatrix}' = \begin{pmatrix} (s_x^{(i)})' \\ (s_y^{(i)})' \\ \sum_{j=1}^{m}\frac{s_x^{(i)}-d_x^{(j)}}{(\|\mathbf{s}^{(i)}-\mathbf{d}^{(j)}\|^{2}+\epsilon)^{\lambda/2}} \\ \sum_{j=1}^{m}\frac{s_y^{(i)}-d_y^{(j)}}{(\|\mathbf{s}^{(i)}-\mathbf{d}^{(j)}\|^{2}+\epsilon)^{\lambda/2}} \end{pmatrix}\quad \forall 1\leq i\leq n,\qquad&\label{eq:sheep_evolution} \\
&\begin{pmatrix} s_x^{(i)}(0) \\ s_y^{(i)}(0) \end{pmatrix} = \mathbf{s}_{0}^{(i)} \quad \forall 1\leq i\leq n,\qquad&
\end{align}

The relative importance of each goal is specified by the hyperparameters $\alpha$ and $\beta$. Note that equation~(\ref{eq:sheep_evolution}) is equivalent to equation (\ref{eq:sheep}) expressed as a first-order ODE; however, the state space is now $4(m+n)$-dimensional. Thus, if we include the co-state, the system is $8(m+n)$-dimensional. Even in the simple case of $m = 2$ and $n=1$, this is a 24-dimensional system. The high number of dimensions increases the spatial complexity of storing a solution and the temporal complexity of calculating the derivative function. Even worse, estimating or directly calculating the Jacobian matrix is signifcantly slower, with the Jacobian matrix containing $24^2 = 576$ entries at each time step.

\begin{figure}
\begin{center}
    \includegraphics[width=0.6\textwidth]{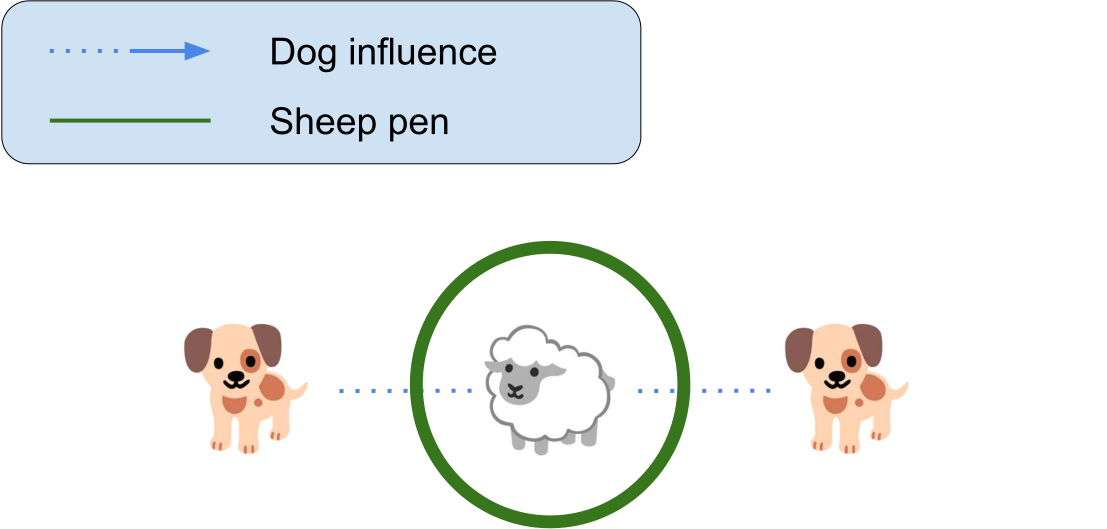}
\end{center}
     \caption[Behavior near $1/\sqrt{5}$]{One potential dog strategy for $m=2$ dogs and $n=1$ sheep. Here, the dogs end by surrounding the sheep at equal distances, forcing the sheep to remain stationary in its pen at the origin.}
\label{fig:dog_strategy} %
\end{figure}
A snapshot of one potential corralling method is given in Figure~\ref{fig:dog_strategy} at time $t=t_f$. Note that our control-based approach does not allow us to impose hard inequality constraints on the state, so we cannot directly enforce that $\|\mathbf{s}^{(i)}(t_f)\|\leq r$ for some final time $t_f$ (specifying the sheep must be contained within a pen of radius $r$ at time $t_f$). Instead, we enforce this as a soft constraint inside our cost functional.

\subsection{Derivation of Optimal Control Using Pontryagin's Maximum Principle}
In solving the optimization problem, we would like to determine the following optimal controls:
\begin{align*}
u_x^{(j)}&\equiv \text{acceleration of dog $j$ in the $x$-direction},\quad\forall 1\leq j\leq m, \\
u_y^{(j)}&\equiv \text{acceleration of dog $j$ in the $y$-direction},\quad\forall 1\leq j\leq m.
\end{align*}
We may do this by calculating the Hamiltonian and deriving the optimal control using Pontryagin's maximum principle. Suppose $\mathbf{p}_d^{(j)}$ denotes the costate corresponding to the position of the $j$th dog and $\mathbf{q}_d^{(j)}$ denotes the costate corresponding to the velocity of the $j$th dog. Suppose we define $\mathbf{p}_s^{(i)}$ and $\mathbf{q}_s^{(i)}$ similarly for sheep. Also, suppose $\mathbf{v}_d^{(j)}$ and $\mathbf{v}_s^{(i)}$ denote the velocity of the $j$th dog and $i$th sheep, respectively. Then the Hamiltonian becomes
\begin{align*}
H &= \sum_{j=1}^m (\mathbf{p}_d^{(j)} \cdot (\mathbf{d}^{(j)})' + \mathbf{q}_d^{(j)} \cdot (\mathbf{v}_d^{(j)})') + \sum_{i=1}^n (\mathbf{p}_s^{(i)} \cdot (\mathbf{s}^{(i)})' + \mathbf{q}_s^{(i)} \cdot (\mathbf{v}_s^{(i)})') &\\
&\quad\quad- \alpha \sum_{i=1}^n\|\mathbf{s}^{(i)}\|_2^2-\beta\sum_{j=1}^m\|\mathbf{d}^{(j)}\|_2^2-\sum_{j=1}^m\|\mathbf{u}^{(j)}\|_2^2\\
&= \sum_{j=1}^m (\mathbf{p}_d^{(j)} \cdot \mathbf{v}_d^{(j)} + \mathbf{q}_d^{(j)} \cdot \mathbf{u}^{(j)}) \\
&\quad\quad+ \sum_{i=1}^n \left(\mathbf{p}_s^{(i)} \cdot \mathbf{v}_s^{(i)} + \mathbf{q}_s^{(i)} \cdot \sum_{j=1}^m \frac{\mathbf{s}^{(i)}-\mathbf{d}^{(j)}}{(\|\mathbf{s}^{(i)}-\mathbf{d}^{(j)}\|_2^2+\epsilon)^{\lambda/2}}\right)\\
&\quad\quad- \alpha \sum_{i=1}^n\|\mathbf{s}^{(i)}\|_2^2-\beta\sum_{j=1}^m\|\mathbf{d}^{(j)}\|_2^2-\sum_{j=1}^m\|\mathbf{u}^{(j)}\|_2^2
\end{align*}
We may now differentiate to get a boundary value problem giving the optimal control. First, differentiating with respect to the control gives
\[
\mathbf{0} = \frac{\partial H}{\partial \mathbf{u}^{(j)}} = \mathbf{q}_d^{(j)} - 2 \mathbf{u}^{(j)}
\]
Therefore,
\[
\mathbf{u}^{(j)} = \tfrac{1}{2} \mathbf{q}_d^{(j)}.
\]
We may also calculate the evolution of the costate:
\begin{align*}
    (\mathbf{p}_d^{(j)})' &=-\frac{\partial H}{\partial \mathbf{d}^{(j)}}\\
    &= \sum_{i=1}^n \mathbf{q}_s^{(i)} \cdot \frac{\partial}{\partial \mathbf{d}^{(j)}} \left(\frac{\mathbf{d}^{(j)} - \mathbf{s}^{(i)}}{(\|\mathbf{d}^{(j)} - \mathbf{s}^{(i)}\|_2^2+\epsilon)^{\lambda/2}} \right) + 2\beta \mathbf{d}^{(j)}\\
    (\mathbf{q}_d^{(j)})' &= -\frac{\partial H}{\partial \mathbf{v}_d^{(j)}}\\
    &= - \mathbf{p}_d^{(j)}\\
    (\mathbf{p}_s^{(i)})' &= -\frac{\partial H}{\partial \mathbf{s}^{(i)}}\\
    &= -\mathbf{q}_s^{(i)} \cdot \sum_{j=1}^m \frac{\partial}{\partial \mathbf{s}^{(i)}} \left(\frac{\mathbf{s}^{(i)} - \mathbf{d}^{(j)}}{(\|\mathbf{s}^{(i)} - \mathbf{d}^{(j)}\|_2^2+\epsilon)^{\lambda/2}} \right) + 2\alpha \mathbf{s}^{(i)}\\
    (\mathbf{q}_s^{(i)})' &= -\frac{\partial H}{\partial \mathbf{v}_s^{(i)}}\\
    &= - \mathbf{q}_s^{(i)}.
\end{align*}
Since the right endpoint is free and there is no final endpoint costs, the costates will zero at the final endpoint. In summary, this gives the BVP:
\begin{align*}
    \begin{pmatrix}
        \mathbf{d}^{(j)} \\
        \mathbf{v}_d^{(j)} \\
        \mathbf{s}^{(i)} \\
        \mathbf{v}_s^{(i)} \\
        \mathbf{p}_d^{(j)} \\
        \mathbf{q}_d^{(j)} \\
        \mathbf{p}_s^{(i)} \\
        \mathbf{q}_s^{(i)}
    \end{pmatrix}' &= \begin{pmatrix}
    \mathbf{v}_d^{(j)} \\
    \frac{1}{2} \mathbf{q}_d^{(j)}\\
    \mathbf{v}_s^{(i)} \\
    \sum_{j=1}^m \frac{\mathbf{s}^{(i)}-\mathbf{d}^{(j)}}{(\|\mathbf{s}^{(i)}-\mathbf{d}^{(j)}\|_2^2+\epsilon)^{\lambda/2}} \\
    \sum_{i=1}^n \mathbf{q}_s^{(i)} \cdot \frac{\partial}{\partial \mathbf{d}^{(j)}} \left(\frac{\mathbf{d}^{(j)} - \mathbf{s}^{(i)}}{(\|\mathbf{d}^{(j)} - \mathbf{s}^{(i)}\|_2^2+\epsilon)^{\lambda/2}} \right) + 2\beta \mathbf{d}^{(j)}\\
     -\mathbf{p}_d^{(j)}\\
     -\mathbf{q}_s^{(i)} \cdot \sum_{j=1}^m \frac{\partial}{\partial \mathbf{s}^{(i)}} \left(\frac{\mathbf{s}^{(i)} - \mathbf{d}^{(j)}}{(\|\mathbf{s}^{(i)} - \mathbf{d}^{(j)}\|_2^2+\epsilon)^{\lambda/2}} \right) + 2\alpha \mathbf{s}^{(i)}\\
     - \mathbf{q}_s^{(i)}
    \end{pmatrix}\\
    \begin{pmatrix}
        \mathbf{d}^{(j)} \\
        \mathbf{v}_d^{(j)} \\
        \mathbf{s}^{(i)} \\
        \mathbf{v}_s^{(i)}
    \end{pmatrix}(0) &= \begin{pmatrix}
        \mathbf{d}^{(j)}_0 \\
        \mathbf{v}_{d,0}^{(j)} \\
        \mathbf{s}^{(i)}_0 \\
        \mathbf{v}_{s,0}^{(i)}
    \end{pmatrix}\\
    \begin{pmatrix}
        \mathbf{p}_d^{(j)} \\
        \mathbf{q}_d^{(j)} \\
        \mathbf{p}_s^{(i)} \\
        \mathbf{q}_s^{(i)}
    \end{pmatrix}(t_f) &= \mathbf{0}\\
    \mathbf{u}^{(j)} &= \tfrac{1}{2} \mathbf{q}_d^{(j)}
\end{align*}
where the indices run over all $1 \leq j \leq m$ and $1 \leq i \leq n$. Finally, note that the above derivative terms can be explicitly calculated using the formula,
\begin{equation} \label{eq:messy_derivative}
\frac{\partial}{\partial \mathbf{x}} \left(\frac{\mathbf{x}}{(\|\mathbf{x}\|^2 + \epsilon)^{\lambda/2}} \right) = (\|\mathbf{x}\|^2 + \epsilon)^{-\lambda/2} (I - \lambda(\|\mathbf{x}\|^2 + \epsilon)^{-1}\mathbf{x}\mathbf{x}^\top).
\end{equation}
\section{Methods and Results}
All the above methods attempt to solve the above boundary value problem. Some methods attempt to find a solution to the exact, nonlinear equation. Other methods attempt to find a solution to an approximation of the problem.

\subsection{Solving Numerically with \texttt{solve\_bvp}}

\subsubsection{Method}

We begin by using the SciPy function \texttt{integrate.solve\_bvp} to numerically solve the original, nonlinear problem. This function is based on a collocation algorithm \cite{kierzenka2001bvp}. Essentially, the solution is approximated with a piecewise cubic polynomial, whereupon the differential equation evaluated on a set of mesh points becomes an algebraic equation that can be solved exactly. Although versatile, the algorithm may require a large number of mesh points, which becomes expensive in terms of both computation and memory. 

Using SymPy, we were also able to verify and find an alternative form of the relatively complicated derivative given in Equation \ref{eq:messy_derivative}.








\subsubsection{Results}

\begin{figure}
\begin{center}
    \includegraphics[width=0.4\textwidth]{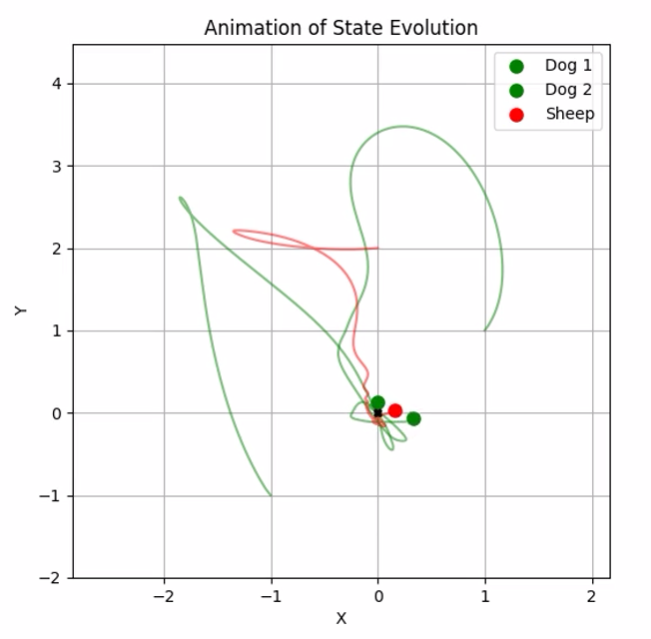}
     \caption[Behavior near $1/\sqrt{5}$]{A visualization of the output from the \texttt{solve\_bvp} method with $m=2$ dogs and $n=1$ sheep. The two dogs position themselves in front of the sheep and progressively steer it toward the origin, circling around it as it nears the destination.}
\label{fig:2Dogs1Sheep}
\end{center}
\end{figure}

Figure \ref{fig:2Dogs1Sheep} displays the results for the  \texttt{solve\_bvp} method. Starting from a random initial guess, the resulting solution has two dogs taking intelligent paths to immediately begin herding the sheep toward the origin. We found that typically the following hyperparameters produced reasonable results: $\lambda = 3, \alpha = 30, \beta = 2, \epsilon = .001$.

Although the visualization appears promising, the solver does not fully converge. This is due to \texttt{solve\_bvp} stopping as soon as the maximum number of mesh nodes is reached, leading the algorithm to return a ``solution" whose derivative poorly approximates the differential equation. Additionally, performance degrades when the numbers of dogs or sheep is increased.

\begin{figure}
\begin{center}
    \includegraphics[width=0.6\textwidth]{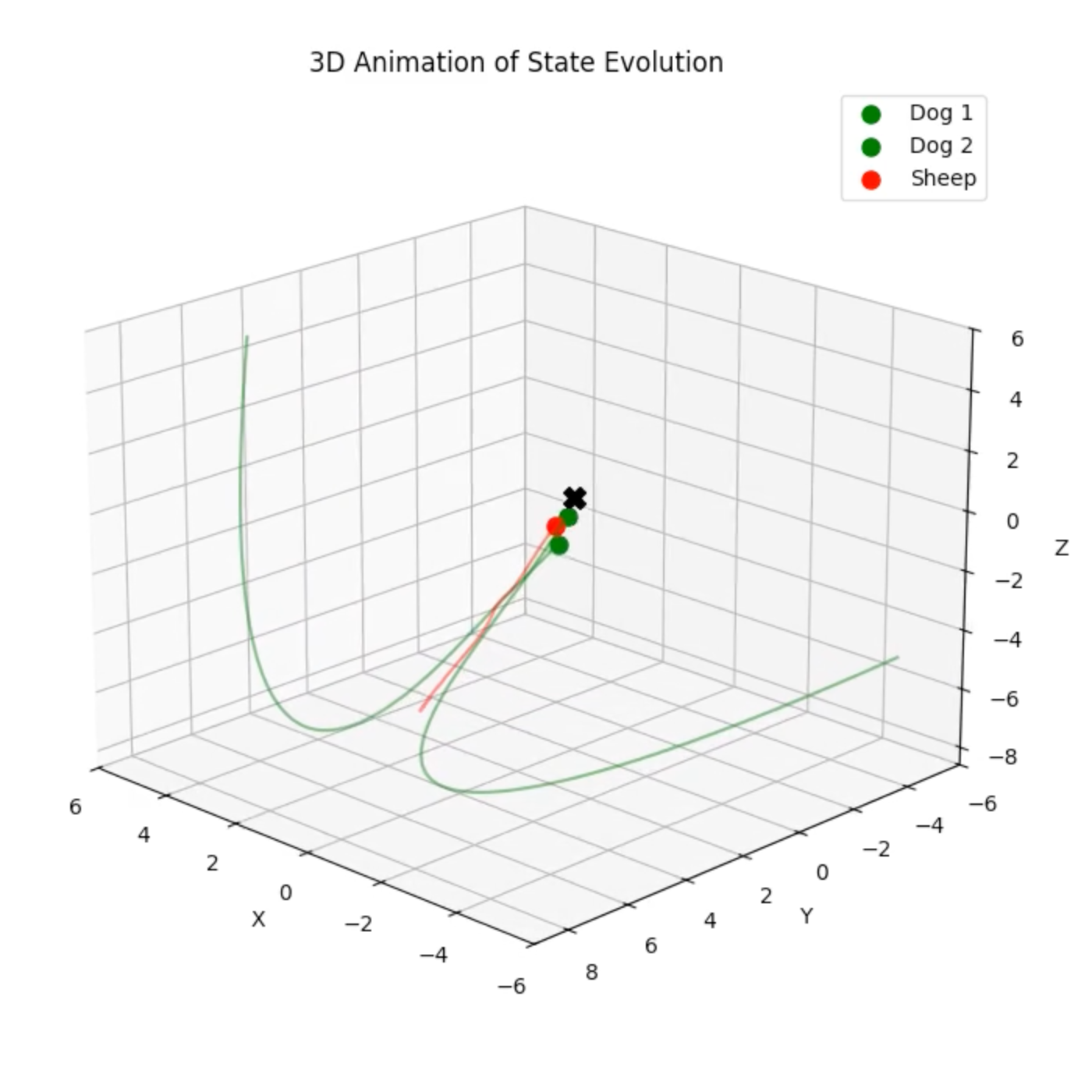}
     \caption[Behavior near $1/\sqrt{5}$]{A visualization of the output from the \texttt{solve\_bvp} method in 3 dimensional space. The two dogs position themselves in behind the sheep and push it toward the origin.}
\label{fig:3d-dog-sheep}
\end{center}
\end{figure}

We can also naturally extend the simulation to three dimensions by augmenting the state equations to include a $z$-dimension term in equations~(\ref{eq:dog_evolution}) and~(\ref{eq:sheep_evolution}). One example result is shown in Figure \ref{fig:3d-dog-sheep}. Though the solution is not interpretable in the context of herding sheep, it could be in the context of robotic control where the herding agents and targets are not constrained to a plane.

\subsection{Solving Numerically with Linearization and
Infinite-Time LQR}

\subsubsection{Method}

One issue we ran into with \texttt{solve\_bvp} was that the solutions typically did not converge and as a result the trajectories did not follow the state space equations. One way to solve this is to change the \texttt{solve\_bvp} to a \texttt{solve\_ivp}. We do this using an iterative Linear Quadratic Regulator.

Notice that the cost functional (\ref{eq:Cost Functional}) is quadratic in all variables, and the dog dynamics (\ref{eq:dog_evolution}) are linear. Although the sheep dynamics (\ref{eq:sheep_evolution}) are nonlinear, we can linearize them around each point along the trajectory, enabling the use of Linear Quadratic Regulator (LQR) techniques to approximate the optimal control at each time step.

To apply LQR effectively, we extend the cost functional to an infinite time horizon. This is important because solving the finite-horizon case requires integrating an initial value problem (IVP) for the $P$ matrix. Since $A(t)=Df(x(t))$ depends on the state and thus varies over time, computing the control would require knowledge of the entire trajectory in advance. By considering the infinite-horizon LQR formulation, we instead solve the Algebraic Riccati Equation (ARE) at each time step, avoiding the need for a full trajectory. Alternatively, one could initialize with a reference trajectory and apply an iterative method, such as iLQR, to update the controls. 

\begin{enumerate} 
    \item Initialize the states of all dogs and sheep. 
    \item Integrate the system using \texttt{solve\_ivp}, and at each time step: 
    \begin{enumerate} 
        \item Linearize the dynamics to compute $A(t)$ around the current state. 
        \item Solve the Algebraic Riccati Equation to obtain the matrix $P$. 
        \item Compute the optimal control using the LQR formulation. 
        \item Evaluate the state derivatives based on the current control and return them to \texttt{solve\_ivp}. 
    \end{enumerate} 
    \item Continue until a predefined final time is reached and return the resulting trajectory. 
\end{enumerate}
We found the following coefficient weights yielded good results:
\begin{table}[h]
\centering
\begin{tabular}{l c}
\hline
Parameter & Value \\ 
\hline
Control effort & 10 \\ 
Sheep's distance from the origin & 10 \\ 
Sheep's velocity & 1 \\ 
Dog's velocity & 0.1 \\ 
Dog's position & 0.2 \\ 
Epsilon & 0.1 \\ 
\hline
\end{tabular}
\caption{Simulation parameters for the Infinite-Time LQR}
\label{Params}
\end{table}

\subsubsection{Results}

\begin{figure}
\begin{center}
    \includegraphics[width=0.6\textwidth]{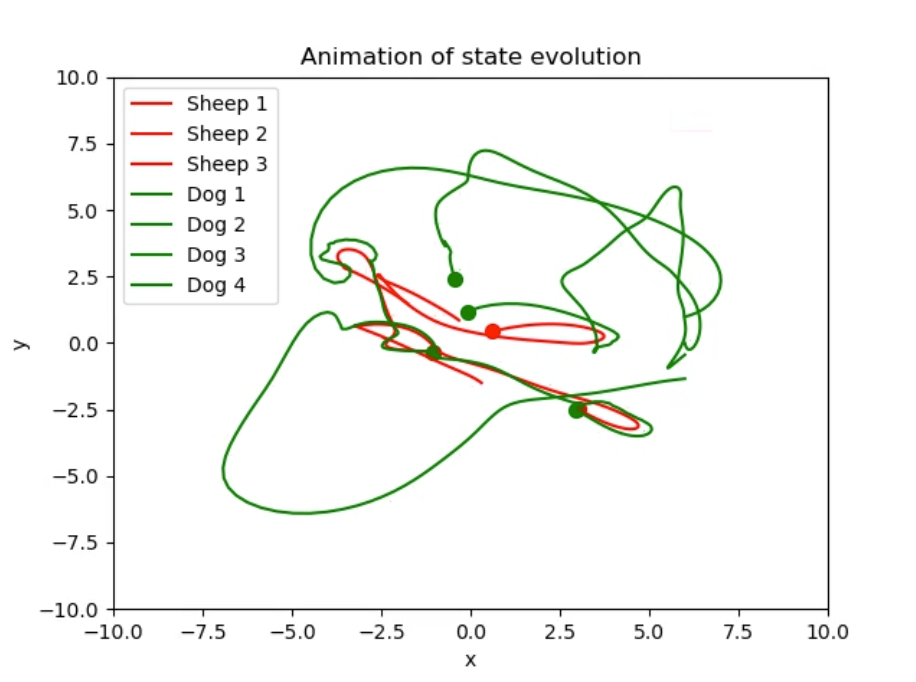}
     \caption[Behavior near $1/\sqrt{5}$]{A visualization of the herding strategy with $m=4$ dogs and $n=3$ sheep. We find that typically one or two dogs will go after each sheep and guide them back to the origin.}
\label{fig:4dog3sheep_linear}
\end{center}
\end{figure}

After implementing the linearized LQR solver, we were able to scale up the simulation to include more agents—successfully handling scenarios with multiple dogs and sheep (see Figure \ref{fig:4dog3sheep_linear}). The dogs were generally able to guide the sheep to the origin, and in many runs, the results were surprisingly effective. Despite the underlying nonlinearity of the system, the LQR-based control was often sufficient to steer the ensemble toward the goal in a coordinated manner. This provided promising empirical evidence that LQR, even in the presence of nonlinearity, could be a viable strategy for multi-agent shepherding tasks.

However, convergence of the solver remained very slow, especially in regions where the sheep were in close proximity to the dogs. This behavior can be traced back to the linearization of the sheep dynamics in equation (\ref{eq:sheep_evolution}). When the sheep are far from the dogs, the influence of any one dog is relatively smooth and the linear approximation is close to the true nonlinear function. But as the sheep approach the dogs, the force terms become sharply nonlinear—often exhibiting near-pole-like behavior—and the linear approximation becomes poor. As a result, the optimal controls computed via LQR may not reflect the actual dynamics accurately, leading to stagnation or oscillatory behavior in the solver. In several cases, dogs would simply ``sit" on top of the sheep, applying ineffective controls while the system evolved only incrementally.

To mitigate the sharp nonlinearities in the interaction terms, we introduced a small $\epsilon$ in the denominator to regularize the near-singular behavior that occurs when dogs approach sheep. While this helped avoid true poles in the dynamics, it did not sufficiently smooth the system to make linearization effective in close-range interactions. As a result, the computed control inputs often failed to reflect the true nonlinear behavior near contact points. This mismatch led to jagged and erratic trajectories—particularly for the dogs—who would sometimes hover near sheep making small, ineffective adjustments. The resulting paths exhibited unnatural oscillations and tight loops, causing the agents to move inefficiently and the solver to progress very slowly.

These issues were further exacerbated when the number of sheep exceeded the number of dogs. With fewer dogs available to manage more sheep, the control authority was spread thin, and the inaccuracies from the linear approximation became even more impactful. In such cases, the dogs frequently got stuck near individual sheep, unable to coordinate effectively as a group to influence the herd. This significantly hindered convergence and revealed a key limitation of the linearized LQR approach: while effective in smoother regions of the state space, it struggles in densely interactive, highly nonlinear scenarios. These observations highlight the need for more robust control strategies capable of handling strong agent-to-agent coupling and sharp local nonlinearities.

\subsection{Improvements to \texttt{solve\_bvp}}
We tried numerous modifications to improve the speed and efficiency of SciPy's \texttt{solve\_bvp}. After much work, we were able to get it to successfully converge. The following ideas were implemented:
\begin{enumerate}
    \item Initially, we attempted to code up a naive implementation of the direct (single) shooting method. Using this method, we choose a random initial costate, and the differential equation is run above as an IVP problem. However, at the same time as the regular differential equation, we run a differential equation on the Jacobian of the derivative (multiplied by a certain matrix). When multiplied by the right matrix, the final value of this Jacobian then gives the Jacobian of the final costate with respect to the initial costate. Using this Jacobian, the initial costate is iteratively updated until a choice is made that gives the desired final costate (zero). Ultimately, our implementation quickly diverged.
    \item We coded the above differential equation in JAX and used JAX to get the Jacobian matrix of the derivative at every time step. This improves \texttt{solve\_bvp}'s performance in both speed and accuracy, as the algorithm no longer needs to estimate the Jacobian (which takes more time and is less accurate).
    \item We increased the maximum number of mesh nodes \texttt{solve\_bvp} could create before terminating. Furthermore, instead of simply terminating the algorithm when the maximum number of mesh nodes was reached, we evaluated the final solution estimate on the original set of mesh nodes and used the result as the initial guess of a new \texttt{solve\_bvp} function call. We then called \texttt{solve\_bvp} somewhere between 10 to 30 times (terminating early if it converged).
    \item We modified the initialization. Namely, we initialized the dog positions randomly on a circle of a given radius centered at the origin (usually a radius of 2 or 10). We also randomly initialized the sheep positions to be on a circle of a smaller radius (usually 0.5 or 1). We also tried randomly initializing the initial velocities as either zero or small Gaussian noise. 
    \item We significantly changed the initial guess. Namely, we had each of the sheep move in a linear line towards the origin and we had the dogs move along curves of the form: 
    \[
    (r_0(1-t/t_f) + r_ft/t_f)e^{i(\theta_0(1-t/t_f)+ \theta_ft/t_f)}
    \] The terms $r_0$ and $\theta_0$ describe the initial position of the dog in polar coordinates. The term $r_f$ was chosen to be a radius notably larger than the initial radius of the sheep (such as 1, 5, or 20). The term $\theta_f$ was chosen to be one of the angles of a corresponding root of unity that in a certain sense was closest to the dog. We also used JAX to calculate the velocity and acceleration of these curves, and initialized the velocity and dog co-velocity accordingly (as the dog co-velocity is proportional to the dog acceleration).
    \item We tried various methods of costate initialization. Beyond setting them as zeros and small Gaussian noise, we also tried running the \texttt{solve\_bvp} algorithm a few times, but changing the state back to the initial state guess each iteration. In other words, we used \texttt{solve\_bvp} to update the costate but not the state. After a few iterations, we then ran \texttt{solve\_bvp} like normal. In the end, this method did not end up being effective.
    \item We modified the $\beta\|\mathbf{d}^{(j)}\|^2$ term in the cost functional. Originally, we changed the distance to be from a given root of unity rather than the origin. However, this incentives every dog to end up at a predetermined root of unity, even if the dog is on the opposite side of the plane. Since crossing the plane often causes the sheep to greatly accelerate, this cost functional ultimately did not work. However, we also tried the cost $\tfrac{1}{2}\beta(\|\mathbf{d}^{(j)}\|^2-1)^2$. This cost more generally incentivizes the dogs to be anywhere on the unit circle. This more flexible cost proved to work much better.
    \item We modified the hyperparameters. Namely, we set $t_f$ to be 2 or 5, as this proved to converge more easily. We also changed $\alpha$ and $\beta$ to $1.$ and $0.02$ respectively, so that the acceleration cost of the dogs was weighted more and the position of the dogs was not nearly as significant. Finally, we initially set $\epsilon$ to be $10^{-3}$, but this often resulted in erratic behavior. Changing $\epsilon$ to $10^{-1}$ instead proved to be much more stable.
\end{enumerate}
In the attached files, we include a few animations of solutions generated by this improved use of \texttt{solve\_bvp}.

\section{Conclusion} 

In this work, we developed state evolution equations and a cost functional to address the age old question: \textit{What are the optimal routes for \( m \) dogs to take to corral \( n \) sheep into their pen?} Our formulation captures the essential dynamics of the system and lays the groundwork for further exploration.

Leveraging Pontryagin's Maximum Principle, we derived the optimal equations of motion for a group of dogs herding a group of sheep. To approximate the solution in the case of two dogs and one sheep, we employed two numerical methods: a collocation method implemented via SciPy’s \texttt{solve\_bvp}, and an iterative LQR approach with an infinite time horizon. While both methods captured the desired long-term behavior, the resulting trajectories often failed to exhibit local optimality. This limitation stems from convergence issues that led to inaccurate approximations of the underlying differential equations.

To mitigate these issues, we implemented several improvements. These include fine-tuning the hyperparameters, tryingda the shooting method, employing automatic differentiation, and strategically choosing the initial guesses. These enhancements enabled our model to successfully converge to a desirable solution.

Nonetheless, considerable room for exploration remains. To further enhance the effectiveness of the strategy, future research could investigate alternative cost functionals—such as those that minimize total herding time or penalize only the sheep’s final positions. Additionally, exploring more advanced or specialized numerical methods may help improve convergence, stability, and scalability in more complex herding scenarios.

\section*{Acknowledgments}

We thank Dr. Emily Evans and the ACME program from Brigham Young University for supporting this research project.

\bibliographystyle{plain}
\bibliography{refs}

\begin{thebibliography}{1}

\bibitem{PhysRevE.104.044613}
P.~Forg\'acs, A.~Lib\'al, C.~Reichhardt, and C.~J.~O. Reichhardt.
\newblock Active matter shepherding and clustering in inhomogeneous environments.
\newblock {\em Phys. Rev. E}, 104:044613, Oct 2021.

\bibitem{kierzenka2001bvp}
Jacek Kierzenka and Lawrence~F Shampine.
\newblock A bvp solver based on residual control and the maltab pse.
\newblock {\em ACM Transactions on Mathematical Software (TOMS)}, 27(3):299--316, 2001.

\bibitem{Lama2024}
Andrea Lama and Mario di~Bernardo.
\newblock Shepherding and herdability in complex multiagent systems.
\newblock {\em Physical Review Research}, 6(3):L032012, 2024.

\bibitem{Strombom2014}
Daniel Strömbom, Richard~P. Mann, Alan~M. Wilson, Stephen Hailes, A.~Jennifer Morton, David J.~T. Sumpter, and Andrew~J. King.
\newblock Solving the shepherding problem: heuristics for herding autonomous, interacting agents.
\newblock {\em Journal of the Royal Society Interface}, 11(100):20140719, 2014.

\end{thebibliography}

\appendix

\newpage
\appendix

\end{document}